\documentclass[11pt]{amsart}
\usepackage{amsmath}
\usepackage{amssymb}
\usepackage{tabularx}
\usepackage{enumerate}
\usepackage{graphicx}
\usepackage{texdraw}

\usepackage{mathrsfs}
\usepackage{amsfonts,amssymb,amsmath}

\makeatletter
\@addtoreset{equation}{section}

\makeatother
\newtheorem{thm}{Theorem}[section]
\newtheorem{lem}[thm]{Lemma}

\begin{document}
\title[Traveling Wave Solutions of Competitive Models]{Traveling Wave Solutions of Competitive Models with
Free Boundaries$^\dag$}
 \thanks{$\dag$ This work is supported in part by NSFC (No. 11271285).}
\author[J. Yang and B. Lou]{Jian Yang$^\ddag$ and Bendong Lou$^\ddag$}
\thanks{$\ddag$ Department of Mathematics, Tongji
University, Shanghai 200092, China.}
\thanks{{\bf Emails:} {\sf yangjian86419@126.com} (J. Yang), {\sf
blou@tongji.edu.cn} (B. Lou) }

\begin{abstract} We study two systems of reaction diffusion equations with monostable or bistable type of nonlinearity
and with free boundaries. These systems are used as multi-species competitive model.
For two-species models, we prove the existence of a traveling wave solution which consists of two semi-waves
intersecting at the free boundary.
For three-species models, we also prove the existence of a traveling wave solution which, however,
 consists of two semi-waves and one compactly supported wave in between, each intersecting with its neighbor at the free boundary.
\end{abstract}

\maketitle

\section{Introduction}

In this paper, we study the following two systems:
\begin{equation}\label{erlianbo}
\left\{
\begin{array}{l}
\phi''+c\phi'+f(\phi)=0,\quad x\in(-\infty,0],\\
\psi''+c\psi'+g(\psi)=0,\quad x\in[0,\infty),\\
\phi(0)=\psi(0)=0,\\
\phi(-\infty)=\psi(+\infty)=1,\\
c=-\alpha \phi'(0)-\beta \psi'(0),
\end{array}
\right.
\end{equation}
and
\begin{equation}\label{sanlianbo}
\left\{
\begin{array}{l}
\phi_{1}''+c\phi_{1}'+f_1(\phi_1)=0, \quad x\in(-\infty,0],\\
\phi_{2}''+c\phi_{2}'+f_2(\phi_2)=0,\quad x\in[0,h],\\
\phi_{3}''+c\phi_{3}'+f_3(\phi_3)=0,\quad x\in[h,\infty),\\
\phi_1(0)=\phi_2(0)=\phi_2(h)=\phi_3(0)=0, \\
\phi_1(-\infty)=\phi_3(+\infty)=1,\\
c=-\alpha \phi_1'(0)-\beta_l \phi_2'(0),\\
c=-\beta_r \phi_2'(h)-\gamma \phi_3'(0),
\end{array}
\right.
\end{equation}
where $\alpha, \beta, \beta_l, \beta_r, \gamma$ are positive constants,
$f, g, f_1, f_2$ and $f_3$ are monostable or bistable types of nonlinearities and
$c$ is a constant to be determined together with the unknowns $\phi, \psi, \phi_1$, etc..
In what follows, we say that $f$ is a monostable type of nonlinearity ($f$ is of (f$_M$) type, for short),
if $f\in C^1([0,\infty))$ and
\begin{equation*}
f(0)=0<f'(0), \quad f(1)=0> f'(1),\quad (1-s)f(s) >0 \ \mbox{for } s>0, s\not= 1;
\end{equation*}
we say that $f$ is a bistable type of nonlinearity ($f$ is of (f$_B$) type, for short),
if
\begin{equation*}
\hskip 10mm
\left\{
 \begin{array}{l}
 f\in C^1([0,\infty)), \  f(0)=0> f'(0), \ f(1)=0 > f'(1),\ \int_0^1 f(s) ds >0,\\
  \ f(\cdot) <0 \mbox{ in } (0,\theta)\cup (1,\infty),\medskip  f(\cdot)>0 \mbox{ in } (\theta,1) \mbox{ for some } \theta \in (0,1).
\end{array}
\right.
\end{equation*}
A typical monostable $f$ is $f(u) =u(1-u)$, and a typical bistable $f$ is
$f(u)=u(u-\theta)(1-u)$ with $\theta \in (0,\frac{1}{2})$.
It is known that the equation \eqref{erlianbo}$_1$ has monotonically decreasing traveling
front on $\mathbb{R}$ when $c=c^*_f$, where $c^*_f>0$ is the minimal traveling speed when $f$ is of {\rm (f$_M$)} type,
or the unique traveling speed when $f$ is of {\rm (f$_B$)} type (c.f. section 2). Similarly, the equation
\eqref{erlianbo}$_2$ has monotonically increasing traveling front when $c=c^*_g $, where $c^*_g <0$ is the
maximal speed when $g$ is of {\rm (f$_M$)} type, or the unique speed when $g$ is of {\rm (f$_B$)} type
(c.f. section 2).

On the problem \eqref{erlianbo} we have the following main result.

\begin{thm}\label{thm:existence1}
Assume that $f$ is of {\rm (f$_M$)} or {\rm (f$_B$)} type, $g$ is of {\rm (f$_M$)} or {\rm (f$_B$)} type.

\begin{enumerate}
\item[\rm (i)] Let $\alpha>0$ be a given constant. Then for any $c\in (c^*_g, \hat{c}_f)$,
where $\hat{c}_f >0$ depends only on $\alpha$ and $f$, there exists a unique $\beta(c)>0$
such that \eqref{erlianbo} has a unique solution $(\phi, \psi, c)$.
Moreover, $\beta(c)$ is continuous and strictly decreasing in $c\in (c^*_g, \hat{c}_f)$ and
\begin{equation}\label{thm1-1}
c\to \hat{c}_f \ \Leftrightarrow \ \beta\to 0,\qquad
c\to c^*_g \ \Leftrightarrow \ \beta\to \infty,\qquad c>0\ \Leftrightarrow \ \beta< \tilde{\beta},
\end{equation}
where $\tilde{\beta} := \alpha (\int_0^1f(s)ds / \int_0^1g(s)ds)^{1/2}$.

\item[\rm (ii)] Let $\beta>0$ be a given constant. Then for any $c\in (\hat{c}_g, c^*_f)$,
where $\hat{c}_g <0$ depends only on $\beta$ and $g$, there exists a unique $\alpha(c)>0$
such that \eqref{erlianbo} has a unique solution $(\phi, \psi, c)$.
Moreover, $\alpha(c)$ is continuous and strictly increasing in $c\in (\hat{c}_g, c^*_f)$ and
\begin{equation}\label{thm1-2}
c\to \hat{c}_g \ \Leftrightarrow \ \alpha\to 0,\qquad
c\to c^*_f \ \Leftrightarrow \ \alpha\to \infty,\qquad c>0\ \Leftrightarrow \ \alpha > \tilde{\alpha},
\end{equation}
where $\tilde{\alpha} := \beta (\int_0^1 g(s)ds / \int_0^1 f(s)ds)^{1/2}$.
\end{enumerate}
\end{thm}

\noindent
This theorem indeed implies that, for any $\alpha, \beta>0$ problem \eqref{erlianbo}
has a unique solution $(\phi(\alpha,\beta)$, $\psi(\alpha,\beta), c(\alpha,\beta))$, and
\eqref{thm1-1} holds when $\alpha$ is fixed,
\eqref{thm1-2} holds when $\beta$ is fixed. This conclusion is an analogue of \cite[Theorem 1.1]{CC}.

On the problem \eqref{sanlianbo} we have the following result.

\begin{thm}\label{thm:existence2}
Assume that $f_1, f_2, f_3$ are of {\rm (f$_M$)} or {\rm (f$_B$)} type.
Let $\alpha, \gamma>0$ be given constants, $\sigma\in(0,1)$ (in case $f_2$ is of {\rm (f$_M$)} type),
or $\sigma\in (\bar{\theta},1)$ (in case $f_2$ is of {\rm (f$_B$)} type) be a given constant.
Then there exist $c_- <0 <c_+$ depending only on $f_1, f_2, f_3, \alpha,\gamma$ and $\sigma$ such that for any
$c\in (c_-, c_+)$, there exists a unique pair $(\beta_l (c), \beta_r(c))$,
$\beta_l(c)$ (resp. $\beta_r(c)$) is continuous and strictly decreasing (resp. increasing) in $c$,
such that problem \eqref{sanlianbo} has
solution $(\phi_1, \phi_2, \phi_3, c)$ with $\|\phi_2\|_{L^\infty} =\sigma$ when $\beta_l=\beta_l(c)$ and $\beta_r=\beta_r(c)$.

Moreover, $c >0$ iff $\beta_l(c) < \tilde{\beta}_l$, or iff $\beta_r (c)> \tilde{\beta}_r$,
where
\begin{equation}\label{tilde beta}
\tilde{\beta}_l := \frac{ \alpha\sqrt{\int_0^1f_1(s)ds}}{ {\sqrt{\int_0^{\sigma}f_2(s)ds}}}, \quad
\tilde{\beta}_r := \frac{ \gamma\sqrt{\int_0^1f_3(s)ds}}{ {\sqrt{\int_0^{\sigma}f_2(s)ds}}}.
\end{equation}
\end{thm}


Problem \eqref{erlianbo} arises in the study of traveling wave solutions of the following system of
reaction diffusion equations:
\begin{equation}\label{p}
\left\{
\begin{array}{l}
u_{t}=u_{xx}+f(u),\quad  \quad \quad\quad\quad\quad \quad x<s(t),\ t>0,\\
v_{t}=v_{xx}+g(v),\quad \quad \quad\quad\quad\quad \quad \ x>s(t),\ t>0,\\
u(x,t)=v(x,t)=0,\quad \quad \quad \quad \quad \ x=s(t),\ t>0,\\
s'(t)=-\alpha u_{x}(x,t)-\beta v_{x}(x,t),\quad x=s(t),\ t>0,\\
s(0)=0,\ u(x,0)=u_{0}(x)(x<0),\ v(x,0)=v_{0}(x)(x>0),
\end{array}
\right.
\end{equation}
where $x=s(t)$ is the free boundary to be determined together with $u$ and $v$, $f,g\in C^1$
satisfying $f(0)=g(0)=0$.
In population ecology, the appearance of regional partition of multi-species through strong competition
is one interesting phenomena. In \cite{MYY85,MYY86,MYY87}, Mimura, Yamada and Yotsutani used problem \eqref{p} to describe regional
partition of two species, which are struggling on a boundary to obtain their own habitats. Among others,
they obtained the global existence, uniqueness, regularity and asymptotic behavior of solutions for the problem.
Later \cite{CDHMN, DHMP, HIMN, MN} studied similar strong competitive models. Recently Du and Lin \cite{DuLin}
and Du and Lou \cite{DuLou} studied a free boundary problem, which is essentially the problem
\eqref{p} in case $v\equiv 0$. They constructed some semi-waves to characterize the spreading of
$u$ which represents the density of a new species.
Motivated by these works, Chang and Chen \cite{CC} recently study the traveling wave solution
of \eqref{p} (i.e. problem \eqref{erlianbo}) with logistic type of nonlinearities:
$$
f(u)= u(1-u),\quad g(v)=v(1-v).
$$
They obtain the existence and uniqueness of traveling wave solution, similar as our Theorem \ref{thm:existence1}
but for logistic type of $f$ and $g$.  One of our purpose in this paper is to study problem
\eqref{p} for general monostable or bistable type of nonlinearity. In what follows, when $f$ and $g$ are of
(f$_M$) type and (f$_B$) type, respectively, we call the solution of \eqref{erlianbo} a {\it MB-type} traveling wave
solution for convenience. {\it MM-type}, {\it BM-type} and {\it BB-type} of traveling wave solutions are defined
similarly (see Figure 1). Thus \cite{CC} presented a special {\it MM-type} traveling wave solution, while our Theorem \ref{thm:existence1}
gives all these four types of traveling wave solutions.

When three (or more) species are involved in contesting the habitats, one
should consider the following competitive model:
\begin{equation}\label{p1}
\left\{
\begin{array}{l}
u_{1t}=u_{1xx}+f_1(u_1),\hskip 42mm x<s_l(t),\ t>0,\\
u_{2t}=u_{2xx}+f_2(u_2),\quad  \quad\quad\quad\quad\quad\quad\quad\quad\quad\quad s_l(t)<x<s_r(t),\ t>0,\\
u_{3t}=u_{3xx}+f_3(u_3),\quad \quad\quad\quad\quad\quad\quad\quad\quad\quad\quad x>s_r(t),\ t>0,\\
u_1(x,t)=u_2(x,t)=u_2(\widetilde{x},t)=u_3(\widetilde{x},t)=0, \quad x=s_l(t),\ \widetilde{x}=s_r(t),\ t>0,\\
{s_l}'(t)=-\alpha u_{1x}(x,t)-\beta_l u_{2x}(x,t),\quad \quad \quad \quad\quad \ x=s_l(t),\ t>0,\\
{s_r}'(t)=-\beta_r u_{2x}(\widetilde{x},t)-\gamma u_{3x}(\widetilde{x},t), \quad \quad \quad \quad\quad  \widetilde{x}=s_r(t),\ t>0,\\
u_2(x,0)=u_{20}(x) (0<x<h), \quad s_l(0)=0,\quad s_r(0)=h (0<h<\infty),\\
u_1(x,0)=u_{10}(x) (x<0),\quad u_3(x,0)=u_{30}(x) (x>h).
\end{array}
\right.
\end{equation}
Our problem \eqref{sanlianbo} is nothing but the problem for the traveling wave solutions of \eqref{p1}:
$$
\begin{array}{c}
u_1(x,t)= \phi_1(x-ct)\ (x\leq ct),\quad u_2(x,t)= \phi_2(x-ct)\ (ct\leq x\leq ct+h),\\
u_3(x,t)= \phi_3(x-ct-h)\ (x\geq ct+h),\quad s_l(t)=ct,\quad s_r(t)=ct+h.
\end{array}
$$
As above, if $f_1, f_2$ and $f_3$ are of (f$_B$), (f$_M$) and (f$_B$) types of nonlinearities, respectively,
we call the solution of \eqref{sanlianbo} a {\it BMB-type} traveling wave solution for convenience.
Similarly, one can define {\it MMM-type}, {\it MBM-type} and other types of traveling wave solutions
(see Figure 1). Our Theorem \ref{thm:existence2} indeed includes all of these types.
We point out that similar conclusions as in Theorems  \ref{thm:existence1} and \ref{thm:existence2}
remain true for the models including four or more species. In other words, for such a model, one can
construct a traveling wave which consists of two semi-waves and several compactly supported waves
in between, each intersecting with its neighbor at the free boundary.

In section 2, we give some basic phase plane analysis and
prove Theorem \ref{thm:existence1}. In section 3 we prove Theorem \ref{thm:existence2}.


\section{The Proof of Theorem 1.1}\label{sec:existenceBBtype}

In this section we prove Theorem \ref{thm:existence1} for {\it BB-type} traveling wave solutions,
that is, for the case where both $f$ and $g$ are of (f$_B$) type. Other
types can be proved similarly.

\subsection{Semi-waves and Phase Plane Analysis}
\label{sec:semi-waves}
As in \cite{DuLou}, we call $\phi(z)$ a semi-wave with speed $c$ if $(c,\phi(z))$ satisfies
\begin{equation}\label{semi-wave}
\left\{
\begin{array}{l}
\phi'' +c\phi' +f(\phi)=0 \quad \mbox{ for } z\in (-\infty,0],\\
\phi(0)=0, \ \phi(-\infty)=1, \ \phi(z)>0\ \mbox{ for } z\in (-\infty,0).
\end{array}
\right.
\end{equation}

The equation in \eqref{semi-wave} can be written in the equivalent form
\begin{equation}
\label{newq-p} \phi'=: \Phi,\quad  \Phi'=-c\Phi-f(\phi).
\end{equation}
As long as $\Phi<0$, $\Phi$ can be regarded as a function of $\phi$ which satisfies
\begin{equation}
\label{newP11} \frac{d\Phi(\phi)}{d\phi}=-c-\frac{f(\phi)}{\Phi}.
\end{equation}
For any $\omega<0$, one can consider this equation with initial data $\Phi(\phi)|_{\phi=0} =\omega$.
By a phase plane analysis (c.f. \cite{AW1, AW2, DuLou}), we see that for each
$\omega <0$, there exists exactly one $c=c(\omega)$ such that the solution of \eqref{newP11} satisfies
$\Phi(\phi)\to 0$ as $\phi\to 1^-$. This solution corresponds to a trajectory of \eqref{newq-p}
through $(0,\omega)$ and $(1,0)$ in the semistrip
\begin{equation*}
S_{\phi}=\{(\phi,\Phi):0<\phi<1,\Phi<0\}
\end{equation*}
in $\phi \Phi$-phase plane. This trajectory gives a unique solution $(c(\omega), \phi(z;c(\omega)))$ for
the problem \eqref{semi-wave} with $\phi'(0;c(\omega)) =\omega$. Moreover, as in \cite{AW1, AW2, DuLou},
$c(\omega)$ is continuous and increasing in $\omega\in (-\infty,0)$ and
\begin{equation}\label{c omega limit}
c(\omega)\to c^*_f \mbox{ as } \omega \to 0,\qquad
c(\omega)\to -\infty \mbox{ as } \omega \to -\infty,
\end{equation}
where $c^*_f >0$ is the unique traveling speed of the following problem
\begin{equation}\label{decreasing TW}
\left\{
\begin{array}{l}
\phi'' +c\phi' +f(\phi)=0 \quad \mbox{ for } z\in \mathbb{R},\\
\phi(-\infty)=1, \ \phi(\infty)=0,\ \phi'(z)<0 \mbox{ for } z\in \mathbb{R}.
\end{array}
\right.
\end{equation}

In summary we have the following result.

\begin{lem}\label{lem:phi Phi}
For any $c\in (-\infty, c^*_f)$, problem \eqref{semi-wave} has a unique solution $(c, \phi (z;c))$.
Moreover, $\phi' (0;c)\ (=\omega)$ is continuous and increasing in $c\in (-\infty,c^*_f)$.
\end{lem}

We also need to consider a similar semi-wave $\psi$ with increasing profile:
\begin{equation}\label{semi-wave2}
\left\{
\begin{array}{l}
\psi'' + c\psi' +g(\psi)=0,\;z\in [0,\infty),\\
\psi(0)=0,\;\psi(\infty)=1,\;\psi(z)>0 \mbox{ for } z\in (0,\infty).
\end{array}\right.
\end{equation}
Denote $c^*_g $ the unique traveling speed of the following problem
\begin{equation}\label{semi-wave-3}
\left\{
\begin{array}{l}
\psi'' +c\psi' +g(\psi)=0 \quad \mbox{ for } z\in \mathbb{R},\\
\psi(-\infty)=0, \ \psi(\infty)=1,\ \psi'(z)>0 \mbox{ for } z\in \mathbb{R}.
\end{array}
\right.
\end{equation}
Then $c^*_g <0$ and in a similar way as above one can obtain the following result.

\begin{lem}\label{lem:psi Psi}
For any $c\in (c^*_g, \infty)$, problem \eqref{semi-wave2} has a unique solution $(c, \psi (z;c))$.
Moreover, $\psi' (0;c)$ is continuous and increasing in $c\in (c^*_g, \infty)$.
\end{lem}

\subsection{Proof of Theorem \ref{thm:existence1}}
We only prove (i). The proof of (ii) is similar.

Now $\alpha>0$ is given. Since $(\alpha \phi'(0;c)+ c)|_{c=0} <0$ and
$(\alpha \phi'(0;c)+ c)|_{c\to c^*_f} >0$, we know by Lemma \ref{lem:phi Phi}
that there exists a unique $\hat{c}_f \in (0,c^*_f)$ such that
\begin{equation}\label{D <0}
\alpha \phi'(0;\hat{c}_f)+ \hat{c}_f =0\quad \mbox{ and } \quad
\alpha \phi'(0;c)+ c <0 \mbox{ for all } c\in (c^*_g, \hat{c}_f).
\end{equation}

For any $\beta\geq 0$, we consider the function
\begin{equation}\label{functionofc}
D(c;\beta):=\alpha \phi'(0;c)+\beta \psi'(0;c)+c,\quad c\in (c^*_g, c^*_f).
\end{equation}
By Lemmas \ref{lem:phi Phi} and \ref{lem:psi Psi}, $D(c;\beta)$ is continuous and
strictly increasing in $c$.
For any $c\in (c^*_g, \hat{c}_f)$, $D(c;0)<0$ by \eqref{D <0} and $D(c;\infty)=\infty$.
Hence there exists a unique $\beta(c)>0$ such that $D(c;\beta(c))=0$, that is,
$$
\alpha \phi'(0;c)+\beta(c) \psi'(0;c)+c \equiv 0,\quad \forall c\in (c^*_g, \hat{c}_f).
$$
By Lemmas \ref{lem:phi Phi} and \ref{lem:psi Psi} again, we have $\beta(c)$ is continuous
and strictly decreasing in $c\in (c^*_g, \hat{c}_f)$.  Moreover,
as $c\to \hat{c}_f$ we have $\beta(c)\to 0$ by \eqref{D <0}.
As $c\to c^*_g $ we have $\alpha \phi'(0;c)+c = \alpha \phi' (0;{c^*_g})+c^*_g <0$ and
$\psi'(0;c) \to 0$, and so $\beta(c)\to \infty$.

When $c=0$, by integration we have
$$
D(0;\tilde{\beta})=-\alpha\sqrt{2\int_0^1 f(s)ds}+\tilde{\beta}\sqrt{2\int_0^1 g(s)ds} =0, \quad \mbox{ where }
\tilde{\beta} := \frac{ \alpha\sqrt{\int_0^1f(s)ds}}{ {\sqrt{\int_0^1 g(s)ds}}},
$$
that is,  $\beta(0)=\tilde{\beta}$. Therefore, $c>0$ if and only if $\beta<\tilde{\beta}$.
This completes the proof of Theorem \ref{thm:existence1}.
\qed


\section{The Proof of Theorem \ref{thm:existence2} }\label{sec:existenceBMBtype}
In this section we prove Theorem \ref{thm:existence2}, only for {\it BMB-type} traveling wave solution.
{\it MMM-type}, {\it MBM-type} and other types of traveling wave solutions are proved similarly.
So in this section, we assume $f_1,f_2,f_3$ are of (f$_B$), (f$_M$) and (f$_B$) types, respectively.

\subsection{Compactly Supported Traveling Wave}\label{sec:traveling-pulse}
For any given $\sigma\in (0,1)$ (when $f_2$ is of (f$_B$) type, we choose $\sigma \in (\bar{\theta},1)$),
we call $\phi_2(z)$ a {\it compactly supported traveling wave} with speed $c$ and with height $\sigma$ if,
for some $h>0$, the pair $(c,\phi_2(z))$ solves
\begin{equation}\label{travelingpulse}
\left\{
\begin{array}{lll}
\phi_2''(z) +c\phi_2'(z) +f_2(\phi_2(z))=0  \mbox{ for } z\in (0,h),\\
\phi_2(0)=\phi_2(h)=0, \ \phi_2'(0)>0,\phi_2'(h)<0,\\
\phi_2(z)>0\ (z\in (0,h)) \mbox{ and } \|\phi_2\|_{C([0,h])} =\sigma.
\end{array}
\right.
\end{equation}
For such a solution, we will see below, there exists a unique $z_0\in (0,h)$ such that $\phi_2$ is strictly
increasing in $(0,z_0)$ and strictly decreasing in $(z_0, h)$, and $\phi_2(z_0)=\sigma$.

To study the existence of solutions of \eqref{travelingpulse}, we use a phase plane analysis as above.
The equation of $\phi_2$ is equivalent to
\begin{equation}\label{newq-p1}
\phi_2'=: \Phi_2,\; \Phi_2'=-c\Phi_2-f_2(\phi_2).
\end{equation}
This system has many trajectories, depending on $c$, passing through $(\sigma,0)$ on the
$\phi_2\Phi_2$-phase plane.
More precisely, for any $\omega_l>0$, there exists a unique $c_l (\omega_l)$ such that the trajectory
$T_l(\omega_l;\sigma)$ of \eqref{newq-p1} with $c=c_l (\omega_l)$ lies in $\{(\phi_2,\Phi_2):0<\phi_2<\sigma,\Phi_2>0\}$
and passes through $(0,\omega_l)$ and $(\sigma, 0)$.
As in the previous section, $c_l (\omega_l)$ is strictly increasing in $\omega_l\in (0,\infty)$, and
$$
c_l(\omega_l)\to c^*_l  \mbox{ as } \omega_l \to 0,\qquad
c_l (\omega_l)\to \infty \mbox{ as } \omega_l \to \infty,
$$
where $c^*_l$ is a constant depending only on $\sigma$ and $f_2$ such that the trajectory $T^*_l(\sigma)$ of \eqref{newq-p1} with $c=c^*_l$
passes through $(0,0)$ and $(\sigma, 0)$.

To study the sign of $c^*_l$, we consider
another trajectory $T_l(\omega_l;1)$ besides $T_l(\omega_l;\sigma)$. Here $T_l(\omega_l;1)$
is a trajectory lying above the $\phi_2$-axis and passing through $(0,\omega_l)$ and $(1,0)$.
Such a trajectory clearly exists for some $c=c_l(\omega_l;1)< c_l (\omega_l)$, and so
$T_l(\omega_l;1)$ is above $T_l(\omega_l;\sigma)$. Taking limit as $\omega_l\to 0$ we have
$$
T_l(\omega_l;\sigma)\rightarrow T^*_l (\sigma),\quad T_l(\omega_l;1)\rightarrow T^*_l(1),\quad
T^*_l(1) \mbox{ is above } T^*_l(\sigma) \mbox{ and }
$$
$$
c^*_l = \lim\limits_{\omega_l \to 0} c_l (\omega_l) \geq c^*_l(1):= \lim\limits_{\omega_l\to 0} c_l(\omega_l;1),
$$
where $T^*_l(1)$ is the trajectory corresponding to the traveling wave (the solution of
\eqref{semi-wave-3} with $g$ replaced by $f_2$) with maximal speed
$c^*_l(1)\ (<0)$ (c.f. \cite{AW1, AW2}).
Trajectory $T^*_l(\sigma)$ gives a function $\Phi_2 (\phi;\sigma)$ satisfying
\begin{equation*}
\frac{d\Phi_2}{d\phi}= - c^*_l -\frac{f_2(\phi)}{\Phi_2}, \quad  \Phi_2 (0)= \Phi_2(\sigma)=0.
\end{equation*}
Integrating over $[0,\sigma]$ we have
$$
-c^*_l \int_0^\sigma\Phi_2 (s;\sigma)ds=\int_0^\sigma f_2(s)ds.
$$
Similarly the function $\Phi_2(\phi;1)$ given by $T^*_l(1)$ satisfies
$$
-c^*_l (1) \int_0^1 \Phi_2 (s;1)ds=\int_0^1 f_2(s)ds.
$$
Since $T_l^\ast(1)$ lies above $T_l^\ast(\sigma)$, we have
$\int_0^\sigma\Phi_2(s;\sigma)ds<\int_0^1\Phi_2(s;1)ds$ and so
\begin{equation}\label{c*l}
c^*_l = - \frac{\int_0^\sigma f_2(s)ds}{\int_0^\sigma\Phi_2(s;\sigma)ds}
   < - \frac{\int_0^\sigma f_2(s)ds}{\int_0^1\Phi_2(s;1)ds}=L_\sigma
   := c^*_l(1) \frac{\int_0^\sigma f_2(s)ds}{\int_0^1 f_2 (s) ds} <0 .
\end{equation}

Similarly, for any $\omega_r<0$, there exists a unique $c_r (\omega_r)$ such that the trajectory $T_r(\omega_r;\sigma)$ of
\eqref{newq-p1} with $c=c_r (\omega_r)$ lies in $\{(\phi_2,\Phi_2):0<\phi_2<\sigma,\Phi_2 <0\}$ and
passes through $(0,\omega_r)$ and $(\sigma, 0)$. The function $c=c_r(\omega_r)$ is strictly increasing,
its range is $(-\infty, c^*_r)$, where $c^*_r $ is a constant such that the trajectory of \eqref{newq-p1} with $c=c^*_r$
passes through $(0,0)$ and $(\sigma, 0)$ in $\{(\phi_2,\Phi_2):0<\phi_2<\sigma,\Phi_2 <0\}$, and
\begin{equation}\label{c*r}
c^*_r > R_\sigma := -c^*_l (1) \frac{\int_0^\sigma f_2(s)ds}{\int_0^1 f_2 (s) ds} >0 .
\end{equation}

In summary we have the following result.

\begin{lem}\label{lem:existence of phi2}
For any $c\in \mathbf{C}:= (c^*_l , c^*_r)$, problem \eqref{travelingpulse} has a
unique solution $\phi_2 (z;c)$ on $[0,h_c]$ for some $h_c >0$. Moreover,
$\phi'_2(0;c)$ and $\phi'_2(h_c ;c)$ is strictly increasing in $c\in \mathbf{C}$.
\end{lem}

\subsection{The Proof of Theorem \ref{thm:existence2}}\label{sec:propertylr}

For the functions $\phi_1$ and $\phi_3$ in problem \eqref{sanlianbo}, as in
Lemmas \ref{lem:phi Phi} and \ref{lem:psi Psi} we have the following result.

\begin{lem}\label{lem:1-3}
For any $c\in (-\infty, c^*_1)$ with $c^*_1 := c^*_{f_1} >0$, problem \eqref{semi-wave} with $f=f_1$ has a unique solution
$(c,\phi_1(z;c))$, and $\phi'_1(0;c)$ is continuous and strictly increasing in $c$;
For any $c\in (c^*_3,\infty)$ with $c^*_3 := c^*_{f_3}<0$, problem \eqref{semi-wave2} with $g=f_3$ has a unique solution
$(c, \phi_3 (z;c))$, and $\phi'_3 (0;c)$ is continuous and increasing in $c\in (c^*_3,\infty)$.
\end{lem}

Now we study the free boundary conditions in \eqref{sanlianbo} (c.f. Figures 2 and 3).
For any $\sigma\in (0,1)$, define
\begin{equation}\label{lc0}
D_l(c;\beta_l) :=\alpha \phi'_1(0;c)+\beta_l \phi'_2 (0;c)+c.
\end{equation}
Now we consider the domain of the variable $c$. Denote $\hat{c}_1$ the unique root of
$$
D_l(c;0) = \alpha \phi'_1(0;c) + c =0,
$$
then $0<\hat{c}_1 <c^*_1$. We have two cases:
$$
{\bf Case\ 1:}\ 0<\hat{c}_1 \leq c^*_r;\qquad {\bf Case\ 2:}\ 0<c^*_r < \hat{c}_1.
$$
The domain of $c$ in \eqref{lc0} is $(c^*_l, c^*_1)\cap (c^*_l, c^*_r)$ in {\bf Case\ 1},
and it is $(c^*_l, c^*_r)$ in {\bf Case\ 2}.

In {\bf Case 1} (see Figure 2), for any $\beta_l\geq 0$,  since
$D_l(\hat{c}_1 ;\beta_l) \geq 0$ and $D_l (c^*_l;\beta_l) <0$,
$D_l(c;\beta_l)=0$ has a unique root $c = C_l (\beta_l)\in (c^*_l , \hat{c}_1)$. Moreover,
$C_l (\beta_l)$ is strictly decreasing in $\beta_l$, $C_l(\beta_l)\rightarrow \hat{c}_1 >0$ as $\beta_l\rightarrow0$ and
$C_l(\beta_l)\rightarrow c^*_l <0 $ as $\beta_l\rightarrow\infty$.

In {\bf Case 2}, $D_l (c^*_r;0)  <0$ and $D_l (c^*_r;\beta_l)  >0$ for large $\beta_l>0$.
Hence there exists a unique $\beta^0_l >0$ such that
$$
D_l (c^*_r; \beta^0_l) =  \alpha \phi'_1(0;c^*_r)+\beta^0_l \phi'_2 (0;c^*_r)+c^*_r =0.
$$
It is easily seen that $D_l(c;\beta_l)=0$ has a root $c=C_l(\beta_l)$ if and only if $\beta_l \geq \beta^0_l$.

In summary we have
\begin{equation}\label{Cl}
 \frac{\partial C_l(\beta_l)} {\partial \beta_l}<0, \mbox{ and }
 \left\{
  \begin{array}{ll}
  C_l(0)=\hat{c}_1 >0,& \ C_l(\infty)\rightarrow c^*_l <0 \mbox{ in {\bf Case\ 1}},\\
  C_l(\beta^0_l)= c^*_r >0, & \ C_l(\infty)\rightarrow c^*_l <0 \mbox{ in {\bf Case\ 2}}.
   \end{array}
 \right.
\end{equation}
This gives the zeros of $D_l(c;\beta_l)$ in \eqref{lc0}.



Similarly, define
\begin{equation}\label{rc0}
D_r(c;\beta_r) :=\beta_r \phi'_2(h_c ; c)+\gamma \phi'_3 (0;c)+c.
\end{equation}
Denote $\hat{c}_3$ the unique root of
$$
D_r(c;0) = \gamma \phi'_3(0;c) +c =0,
$$
then $c^*_3 < \hat{c}_3  <0$. On their relations with $c^*_l$ we have two cases:
$$
{\bf Case\ I:}\ c^*_l \leq  \hat{c}_3 < 0;\qquad {\bf Case\ II:}\  \hat{c}_3 < c^*_l <0.
$$
So the domain of $c$ in \eqref{rc0} is $(c^*_3, c^*_r)\cap (c^*_l, c^*_r)$ in {\bf Case\ I},
and it is $(c^*_l, c^*_r)$ in {\bf Case\ II} (see Figure 2). In a similar way as above we see that
the unique root $c = C_r (\beta_r)$ of $D_r(c;\beta_r)=0$ satisfies
\begin{equation}\label{Cr}
 \frac{\partial C_r(\beta_r)} {\partial \beta_r}>0, \mbox{ and }
 \left\{
  \begin{array}{ll}
  C_r(0)=\hat{c}_3 <0,& \ C_r(\infty)\rightarrow c^*_r >0 \mbox{ in {\bf Case\ I}},\\
  C_r(\beta^0_r)= c^*_l <0, &  \ C_r(\infty)\rightarrow c^*_r >0 \mbox{ in {\bf Case\ II}},
   \end{array}
 \right.
\end{equation}
where $\beta^0_r >0$ is the unique root of
$$
D_r(c;\beta_r) :=\beta_r \phi'_2(h_{c^*_l} ; c^*_l)+\gamma \phi'_3 (0;c^*_l)+c^*_l =0.
$$
This gives the zeros of $D_r(c;\beta_r)$ in \eqref{rc0}.

\medskip

Now we try to find some $c$ such that both $D_l(c;\beta_l)=0$ and $D_r(c;\beta_r)=0$ hold at
the same time for suitably chosen pair: $(\beta_l, \beta_r)$. This will give a solution of \eqref{sanlianbo}.

First we consider {\bf Case 1} and {\bf Case II} (see Figure 3).
In this case the range of $C_l(\beta_l)$ is $(c^*_l, \hat{c}_1]$, and the range of $C_r(\beta_r)$
is $(c^*_l, c^*_r)$. Therefore,
for any $c\in (c^*_l, \hat{c}_1)$, there exist a unique $\beta_l = \beta_l (c)$ and
a unique $\beta_r=\beta_r (c)$ such that
$$
D_l(c;\beta_l (c)) = D_r(c;\beta_r (c))=0\quad \mbox{ and so } \ C_l(\beta_l (c))=c=C_r(\beta_r(c)).
$$
Moreover, by the definition of $D_l(c;\beta_l)$ (resp. $D_r(c;\beta_r)$) we see that
$\beta_l(c)$ (resp. $\beta_r(c)$) is strictly decreasing (resp. increasing) function of $c$.
Other cases can be studied similarly. In summary we have the following result.

\begin{thm}\label{thm:exist}
In {\bf Case 1} and {\bf Case I} (resp. in {\bf Case 1} and {\bf Case II},
{\bf Case 2} and {\bf Case I}, {\bf Case 2} and {\bf Case II}), for any $c\in (\hat{c}_3, \hat{c}_1)$ (resp.
$c\in (c^*_l, \hat{c}_1)$, $c\in (\hat{c}_3, c^*_r)$, $c\in (c^*_l, c^*_r)$),
there exists a unique positive pair $(\beta_l (c), \beta_r(c))$ such that \eqref{sanlianbo} has a solution.
\end{thm}

Finally we study the sign of $c$. Again we only consider {\bf Case 1} and {\bf Case II}, since
other cases are treated similarly.
 For $\tilde{\beta}_l$ and $\tilde{\beta}_r$ defined in \eqref{tilde beta}, since
$$
D_l(0;\tilde{\beta}_l)=\alpha \phi'_1(0;0)+\tilde{\beta}_l \phi'_2(0;0)=
-\alpha\sqrt{2\int_0^1 f_1(s)ds}+ \tilde{\beta}_l\sqrt{2\int_0^\sigma f_2(s)ds} =0,
$$
we have $\beta_l (0) = \tilde{\beta}_l$. Similarly, $\beta_r (0) = \tilde{\beta}_r$.
Observing Figure 3 one can find that
$$
c>0 \  \Leftrightarrow \  0< \beta_l<\tilde{\beta}_l\  \Leftrightarrow \ \beta_r > \tilde{\beta}_r .
$$
This completes the proof of Theorem \ref{thm:existence2}. \qed

\medskip

{\bf Acknowledgement}. The authors would like to thank Professor C.C. Chen for sending
the paper \cite{CC} to them, thank Professors M. Nagayama,
K.I. Nakamura and Y. Yamada for helpful discussions.


\begin{thebibliography}{123456}

\bibitem{AW1}D.G.~Aronson and H.F. Weinberger,
{\em Nonlinear diffusion in population genetics, conbustion, and
nerve pulse propagation}, in Partial Differential Equations and
Related Topics, Lecture Notes in Math. 446, Springer, Berlin, 1975,
pp. 5--49.

\bibitem{AW2}D.G.~Aronson and H.F. Weinberger,
{\em Multidimensional nonlinear diffusion arising in population
genetics}, Adv. in Math., {\bf 30} (1978), 33--76.

\bibitem{CC} C.H. Chang and C.C. Chen, {\em Traveling wave solution of a
free boundary problem for a two-species competitive model}, Commun. Pure Appl. Anal.,
{\bf 12} (2013), 1065--1074.

\bibitem{CDHMN}E.C.M. Crooks, E.N. Dancer, D. Hilhorst, M. Mimura and H. Ninomiya,
{\em Spatial segregation limit of a competition-diffusion system with
Dirichlet boundary conditions},  Nonl. Anal. RWA, {\bf 5} (2004),  645--665.


\bibitem{DHMP} E.N. Dancer, D. Hilhorst, M. Mimura and L.A. Peletier,
{\em Spatial segregation limit of a competition-diffusion system},
Euro. J. Appl. Math., {\bf 10} (1999), 97--115.


\bibitem{DuLin}Y.~Du and Z.G.~Lin,
{\em Spreading-vanishing dichtomy in the diffusive logistic model
with a free boundary}, SIAM J. Math. Anal., 42 (2010), 377-405.

\bibitem{DuLou} Y. Du and B.D. Lou,
{\em Spreading and vanishing in nonlinear diffusion problems with free boundaries},
J. Eur. Math. Soc., ({\it to appear}), (http://arxiv.org/pdf/1301.5373.pdf).

\bibitem{HIMN}D. Hilhorst, M. Iida, M. Mimura and H. Ninomiya,
{\em A competition-diffusion system approximation to the classical two-phase stefan problem},
Japan J. Indust. Appl. Math., {\bf 18} (2001), 161--180.

\bibitem{MYY85} M. Mimura, Y. Yamada and S. Yotsutani, {\em A free boundary problem in ecology},
Japan J. Appl. Math., {\bf 2} (1985), 151--186.

\bibitem{MYY86} M. Mimura, Y. Yamada and S. Yotsutani, {\em Stability analysis for free boundary problems
in ecology}, Hiroshima Math. J., {\bf 16} (1986), 477--498.

\bibitem{MYY87} M. Mimura, Y. Yamada and S. Yotsutani, {\em Free boundary problems
for some reaction-diffusion equations}, Hiroshima Math. J., {\bf 17} (1987), 241--280.

\bibitem{MN} H. Murakawa and H. Ninomiya,
{\em Fast reaction limit of a three-component reaction.diffusion system},
J. Math. Anal. Appl., {\bf 379} (2011), 150--170.

\end{thebibliography}
\end{document}